\documentclass[12pt]{article} %platex 
\usepackage[final]{graphicx}
\usepackage[dvips]{color}
\usepackage{graphicx}
\usepackage{amsmath}
\usepackage{amssymb}
\usepackage{ascmac}
\usepackage{amsthm}
\usepackage{bm}
\def\Bbb{\bf} %%%

\newcommand\Z{{\Bbb Z}}
\newcommand\Q{{\Bbb Q}}

\newtheorem{lmm}{Lemma}[section]
\newtheorem{thm}[lmm]{Theorem}
\newtheorem{qst}{Question}

\newtheorem{rmk}[lmm]{Remark}
\newtheorem{prp}[lmm]{Proposition}

\def\comment#1{ }

\makeatletter
    
    \@addtoreset{equation}{section}
  \makeatother
\allowdisplaybreaks
\begin{document}
\title{Some special values of hypergeometric series related to central values of automorphic $L$-functions}
\author{Akihito Ebisu
%\thanks{}
}

\maketitle
\begin{abstract}
In a recent work, H.Narita presented problems concerning 
the strict positivity of central values of certain automorphic $L$-functions
in the form of questions regarding special values of the hypergeometric series.
In this paper, we present partial answers to these questions 
using the theory of orthogonal polynomials and three term relations of the hypergeometric series.

Key Words and Phrases: hypergeometric series, three term relation, central value of $L$-function.

2010 Mathematics Subject Classification Numbers: 33C05. 
\end{abstract}
\tableofcontents
\section{Introduction}
In [NOS], 
H.Narita presented problems regarding  
the strict positivity of central values of certain automorphic $L$-functions
in the form of questions regarding special values of the hypergeometric series.
These questions are the following:
\begin{qst}
Is there a pair $(\kappa _1, \kappa _2)\in (4\Z _{>0})^2$ with $\kappa _1 \leq \kappa _2$ satisfying
\begin{gather}
\binom{(\kappa _2+\kappa _1)/2}{(\kappa _2-\kappa _1)/2}{}_2F_1
\left(\frac{-\kappa _2+\kappa _1}{2}, \frac{-\kappa _2+\kappa _1}{2}; \kappa _1 +1;-1\right)
=(-1)^{\kappa _1/4+1}2^{\kappa _2/2 -2}
\end{gather}
other than $(\kappa _1, \kappa _2)=(4,4), (4,8), (4,16)$ and $(12,16)$?
\end{qst}
\begin{qst}
Is there a pair $(\kappa _1, \kappa _2)\in (4\Z _{>0})^2$ with $\kappa _1 \leq \kappa _2$ satisfying
\begin{gather}
\binom{(\kappa _2+\kappa _1)/2}{(\kappa _2-\kappa _1)/2}{}_2F_1
\left(\frac{-\kappa _2+\kappa _1}{2}, \frac{-\kappa _2+\kappa _1}{2}; \kappa _1 +1;-3\right)
=-2^{\kappa _2 -1}?
\end{gather}
\end{qst}
Here, $_{2}F_{1}(a,b;c;x)$ is the hypergeometric series, and it is defined as
\begin{gather*}
_{2}F_{1}(a,b;c;x):=\sum _{n=0}^{\infty}\frac{(a,n)(b,n)}{(c,n)(1,n)}x^n,
\end{gather*}
where $(a,n):=\Gamma (a+n)/\Gamma (a)$.
With regard to Question 1, if a pair $(\kappa _1, \kappa _2)$ does not satisfy the equality (1.1),
then we can say that the central value of 
the automorphic $L$-function corresponding to that pair is strictly positive.
The same holds for Question 2 (cf. Theorem 0.4 in [NOS]).

Let us define the following quantities:
\begin{align*}
\ell _1(n_1,n_2)&:=\binom{2n_2+2n_1}{2n_2-2n_1}{}_2F_1
\left(-2n_2+2n_1, -2n_2+2n_1; 4n_1+1;-1\right),\\
r_1 (n_1,n_2)&:=(-1)^{n_1+1}2^{2n_2 -2},\\
\ell _3(n_1,n_2)&:=\binom{2n_2+2n_1}{2n_2-2n_1}{}_2F_1
\left(-2n_2+2n_1, -2n_2+2n_1; 4n_1+1;-3\right),\\
r_3(n_1,n_2)&:=-2^{4n_2 -1}.
\end{align*}
Then, we rewrite Quesiton 1 and Question 2 as follows:
\begin{qst}
Is there a pair $(n_1, n_2)\in (\Z _{>0})^2$ with $n_1 \leq n_2$ satisfying 
\begin{gather}
\ell _1(n_1,n_2)=r_1(n_1,n_2)
\end{gather}
other than $(n_1, n_2)=(1,1), (1,2), (1,4)$ and $(3,4)$?
\end{qst}
\begin{qst}
Is there a pair $(n_1, n_2)\in (\Z _{>0})^2$ with $n_1 \leq n_2$ satisfying 
\begin{gather}
\ell _3(n_1,n_2)=r_3(n_1,n_2)?
\end{gather}
\end{qst}
In this paper, we consider Questions $3$ and $4$.
In \S 2, we rewrite $\ell_1(n_1,n_2)$ and $\ell_3(n_1,n_2)$ in simple forms. 
In \S 3, we show that there is no pair $(n_1,n_2)$ with $n_1 \leq n_2$
satisfying (1.4) by using the theory of orthogonal polynomials.
Thus, we answer Question 4 (and 2) completely.
In \S 4, we give the following partial answer to Question 3 (and hence 1):
There is no pair $(n_1,n_2)$ with $4n_1 \leq n_2$ satisfying (1.3) except $(n_1,n_2)=(1,4)$.
This result is obtained by using three term relations of the hypergeometric series.

\begin{rmk}
Computer experiments reveal that no $(n_1,n_2)$ with $0<n_1 \leq n_2 <10^3$ satisfies (1.3) 
other than $(n_1, n_2)=(1,1), (1,2), (1,4)$ and $(3,4)$.
Thus, we conjecture that the answer to Question 3 (and 1) is `No'.
\end{rmk}
\section{Rewriting of $\ell_1(n_1,n_2)$ and $\ell_3(n_1,n_2)$}
In this section, we rewrite $\ell_1(n_1,n_2)$ and $\ell_3(n_1,n_2)$ 
in forms that are more convenient to analyze.

The hypergeometric differential equation $E(a,b,c)$ is given by 
$L(a,b,c)y=0$, where 
\begin{gather*}
L(a,b,c):=\partial ^2+\frac{c-(a+b+1)x}{x(1-x)}\partial -\frac{ab}{x(1-x)},
\end{gather*}
$\partial:=d/dx$, and $a, b$ and $c$ are complex variables. 
The equation $E(a,b,c)$ possesses 
\begin{gather*}
_2F_1(a,b;c;x),\ (-x)^{-a}{}_2F_1(a, a+1-c; a+1-b; 1/x)
\end{gather*}
as a solution (cf. 2.9(1) and (9) in [Erd]).
In particular, 
\begin{align*}
&_2F_1(-2n_2+2n_1,-2n_2+2n_1;4n_1+1; x),\\
&(-x)^{2n_2-2n_1}{}_2F_1(-2n_2+2n_1, -2n_2-2n_1; 1; 1/x)
\end{align*}
satisfies $E(-2n_2+2n_1,-2n_2+2n_1,4n_1+1)$.
Because these solutions are polynomials in $x$,
they are equal up to a multiplicative factor.
By equating the respective coefficients of $x^{2n_2-2n_1}$ of these polynomial solutions,
we are able to calculate this factor explicitly. 
We thereby obtain
\begin{multline}
_2F_1(-2n_2+2n_1,-2n_2+2n_1;4n_1+1; x)\\
=\binom{2n_2+2n_1}{2n_2-2n_1}^{-1}(-x)^{2n_2-2n_1}{}_2F_1(-2n_2+2n_1, -2n_2-2n_1; 1; 1/x).
\end{multline}   
Substituting $x=-1$ and $-3$ into (2.1), we have
\begin{align}
\ell _1 (n_1, n_2)&={}_2F_1(-2n_2+2n_1,-2n_2-2n_1;1; -1),\\
\ell _3 (n_1, n_2)&=3^{2n_2-2n_1}{}_2F_1(-2n_2+2n_1,-2n_2-2n_1;1; -1/3).
\end{align}
Further, applying the Pfaff transformation 
\begin{gather*}
_{2}F_{1}(a,b;c;x)=(1-x)^{-a}{}_{2}F_{1}(a,c-b;c;x/(x-1))
\end{gather*}
(cf. 2.9(1) and (3) in [Erd] or Theorem 2.2.5 in [AAR]) to the right hand sides of (2.2) and (2.3), we obtain
\begin{align}
\ell _1 (n_1, n_2)&=2^{2n_2-2n_1}{}_2F_1(-2n_2+2n_1,2n_2+2n_1+1;1; 1/2),\\
\ell _3 (n_1, n_2)&=4^{2n_2-2n_1}{}_2F_1(-2n_2+2n_1,2n_2+2n_1+1;1; 1/4).
\end{align}
These are represented in terms of the Jacobi polynomial $P_{n}^{(a,b)}(x)$ as
\begin{align}
\ell _1 (n_1, n_2)&=2^{2n_2-2n_1}P_{2n_2-2n_1}^{(0,4n_1)}(0),\\
\ell _3 (n_1, n_2)&=4^{2n_2-2n_1}P_{2n_2-2n_1}^{(0,4n_1)}(1/2),
\end{align}
where
\begin{gather*}
P_{n}^{(a,b)}(x):=\frac{(a+1,n)}{(1,n)}{}_2F_1(-n, n+a+b+1; a+1; (1-x)/2).
\end{gather*}
\begin{rmk}
The expansions 
\begin{align*}
\ell _1(n_1,n_2)&=\sum _{i=0}^{2n_2-2n_1}\binom{2n_2-2n_1}{i}\binom{2n_2+2n_1}{i}(-1)^i\\
\ell _3(n_1,n_2)&=3^{2n_2-2n_1}{}\sum _{i=0}^{2n_2-2n_1}\binom{2n_2-2n_1}{i}\binom{2n_2+2n_1}{i}(-1/3)^i
\end{align*}
of the hypergeometric series given in (2.2) and (2.3) show that these are integers.
\end{rmk}

\section{Answer to Question 4}
In this section, we answer Question 4
by applying the following theorem to $(2.7)$:
\begin{thm}(cf. Theorem 7.2 in [Sz]).
Let $w(x)$ be a weight function that is non-decreasing in the interval $[a,b]$, with $b$ finite.
If $\{ p_{n}(x) \}$ is the set consisting of the corresponding orthogonal polynomials, 
then the functions $\{ w(x)\} ^{1/2}|p_n(x)|$ realize their maxima in $[a,b]$ at $x=b$;
that is, 
\begin{gather*}
\{ w(x)\} ^{1/2}|p_n(x)| \leq \{ w(b)\} ^{1/2}|p_n(b)|
\end{gather*}
for any point $x\in [a,b]$.
\end{thm} 
Because the weight function of $P_{n}^{(0,b)}(x)$ is $w(x)=(1+x)^b$,
we are able to apply the above theorem to $(2.7)$, and we thereby deduce
\begin{align*}
|\ell_3(n_1,n_2)|&=4^{2n_2-2n_1}{}|P_{2n_2-2n_1}^{(0,4n_1)}(1/2)|
\leq 4^{2n_2-2n_1}{}\frac{w(1)^{1/2}}{w(1/2)^{1/2}}|P_{2n_2-2n_1}^{(0,4n_1)}(1)|\\
&\leq 4^{2n_2-2n_1}{}\frac{2^{2n_1}}{(3/2)^{2n_1}}\cdot 1
=2^{4n_2}3^{-2n_1}
<2^{4n_2-1}=|r_3(n_1,n_2)|.
\end{align*}
Thus, we obtain a negative answer to Question 4, and hence to Question 2:
\begin{prp}
There is no pair $(\kappa _1, \kappa _2)\in (4\Z _{>0})^2$ with $\kappa _1 \leq \kappa _2$ satisfying (1.2).
\end{prp}
\begin{rmk}
Applying Theorem 3.1 to (2.6), we find $|\ell_1(n_1,n_2)| \leq 2^{2n_2}$.
However, because $|r_1(n_1,n_2)|=2^{2n_2-2}$,
we cannot determine whether $|r_1(n_1,n_2)|$ is larger than $|\ell_1(n_1,n_2)|$.
Indeed, the relation $|\ell_1(n_1,n_2)|=|r_1(n_1,n_2)|$ holds for the four pairs appearing in Remark 1.1,
and $|\ell_1(n_1,n_2)|>|r_1(n_1,n_2)|$ holds for
\begin{gather*}
(n_1,n_2)=(1,3), (2,3), (3,6), (4,6), (6,9), (8,12).
\end{gather*}
However, we also find that $|\ell_1(n_1,n_2)|<|r_1(n_1,n_2)|$ holds
for any $(n_1,n_2)$ satisfying $0<n_1\leq n_2<10^3$ other than the above ten pairs.
Thus, we conjecture that $|\ell_1(n_1,n_2)|<|r_1(n_1,n_2)|$ for all $(n_1,n_2)$ with $13 \leq n_2$.
\end{rmk}

\section{Partial answer to Question 3}
Because we have
$$
_{2}F_1(a,b;c;x)={}_2F_1(b,a;c;x), 
$$ 
(2.2) leads to
\begin{gather}
\ell_1(n_1,n_2)={}_2F_1(-2n_2-2n_1,-2n_2+2n_1;1;-1).
\end{gather}
Also, recall that
\begin{gather}
_2F_1(a, b; a+1-b; -1)=2^{-a}\frac{\Gamma (a+1-b)\Gamma (1/2)}{\Gamma (a/2+1-b)\Gamma (a/2+1/2)}
\end{gather} 
(cf. 2.8(47) in [Erd]).
Thus, using (4.1) and (4.2), $\ell _1(0,n_2)$ can be expressed in the closed form
\begin{gather}
\ell_1(0,n_2)={}_{2}F_1(-2n_2,-2n_2;1;-1)=\frac{(-1)^{n_2}\left\{1\cdot 3\cdot 5\cdots (2n_2-1) \right\}
2^{n_2}}{(1,n_2)}.
\end{gather}
Next, we investigate the forms of $\ell _1(1,n_2)$, $\ell _1(2,n_2)$, $\cdots$, $\ell _1(n_2,n_2)$.
For this purpose, three term relations of the hyprergeometric series play an important role.

\subsection{Properties of three term relations of the hypergeomtric series}
In this subsection, we state properties of three term relations of the hypergeometric series. 

It is well known (cf. \S 1 in [Eb]) that 
for any triplet of integers $(k,l,m)\in\Z^3$,
there exists a unique pair of rational functions 
$(Q(x), R(x)) \in (\Q (a,b,c,x))^2,$ 
where $\Q (a,b,c,x)$ is the field generated over $\Q$ by $a$, $b$, $c$ and $x$,
satisfying
\begin{gather*}
_{2}F_{1}(a+k,b+l;c+m;x)=Q(x){}_{2}F_{1}(a+1,b+1;c+1;x)+R(x){}_{2}F_{1}(a,b,c;x).
\end{gather*}
This relation is called the `three term relation of the hypergeometric series'. 
In the next paragraph, 
we express $Q(x)$ and $R(x)$ in terms of contiguity operators.

First, let us define $\vartheta := x\partial$. Then, we apply the differential operators 
\begin{align*}
&H_1(a,b,c):=\vartheta +a,\, B_1(a,b,c):=\frac{1}{(1-a)(c-a)}(-x(1-x)\partial +(bx+a-c)),\\
&H_2(a,b,c):=\vartheta +b,\, B_2(a,b,c):=\frac{1}{(1-b)(c-b)}(-x(1-x)\partial +(ax+b-c)),\\
&H_3(a,b,c):=\frac{1}{(c-a)(c-b)}((1-x)\partial +c-a-b),\, B_3(a,b,c):=\vartheta +c-1 
\end{align*}
to the hypergeometric series $_2F_1(a,b;c;x)$. This yields
\begin{align*}
H_1(a,b,c){}_2F_1(a,b;c;x)&=a\,{}_2F_1(a+1,b;c;x),\\
B_1(a,b,c){}_2F_1(a,b;c;x)&=1/(a-1)\cdot {}_2F_1(a-1,b;c;x),\\
H_2(a,b,c){}_2F_1(a,b;c;x)&=b\,{}_2F_1(a,b+1;c;x),\\
B_2(a,b,c){}_2F_1(a,b;c;x)&=1/(b-1)\cdot {}_2F_1(a,b-1;c;x),\\
H_3(a,b,c){}_2F_1(a,b;c;x)&=1/c\cdot{}_2F_1(a,b;c+1;x),\\
B_3(a,b,c){}_2F_1(a,b;c;x)&=(c-1){}_2F_1(a,b;c-1;x)
\end{align*}
(cf. Theorems 2.1.1 and 2.1.3 in [IKSY]).
For this reason, we call these operators `contiguity operators'.
Let $H(k,l,m)$ be a composition of these contiguity operators such that
\begin{gather*}
H(k,l,m)_2F_1(a,b;c;x)=\frac{(a,k)(b,l)}{(c,m)}{}_2F_1(a+k,b+l;c+m;x).
\end{gather*}
Further, we define $q(x)$ and $r(x)$ by
\begin{gather}
Q(x)=\frac{a b (c,m)}{c(a,k)(b,l)}q(x),\  R(x)=\frac{(c,m)}{(a,k)(b,l)}r(x).\ 
\end{gather}
Then, $H(k,l,m)$ can be expressed as
\begin{gather} 
H(k,l,m)=p(\partial)L(a,b,c)+q(x)\partial +r(x),
\end{gather}
where $p(\partial)$ is an element of the ring of the differential operators in $x$ over $\Q (a,b,c,x)$
(cf. (2.4) in [Eb]).

In this paragraph, we investigate the properties of $q(x)$ and $r(x)$.
Without loss of generality, we assume that $k\leq l$.
Then, $q(x)$ and $r(x)$ can be expressed as
\begin{align*}
q(x)&=x^{v_0}(1-x)^{v_1}q_0(x),
& q_0(x):{\text{a polynomial of degree $g$ and }} q_0(0)q_0(1)\neq 0,\\
r(x)&=x^{w_0}(1-x)^{w_1}r_0(x),
& r_0(x):{\text{a polynomial of degree $h$ and }} r_0(0)r_0(1)\neq 0,
\end{align*}
where 
\begin{eqnarray*}
  (v_0,v_1,g)= \left\{ \begin{array}{cl}
   (1-m, m+1-k-l, l-1)  & {\text{if }} (k,l,m)\in {\rm (i)}, \\
   (1-m, 1, m-k-1)  & {\text{if }} (k,l,m)\in {\rm (ii)}, \\
   (1, m+1-k-l, l-m-1)  & {\text{if }} (k,l,m)\in {\rm (iii)},  \\
   (1, 1, -k-1)  & {\text{if }} (k,l,m) \in {\rm (iv)},
  \end{array} \right. \\
  (w_0,w_1,h)= \left\{ \begin{array}{cl}
   (1-m, m+1-k-l, l-2)  & {\text{if }} (k,l,m)\in {\rm (i)}, \\
   (1-m, 0, m-k-1)  & {\text{if }} (k,l,m)\in {\rm {(ii)}}, \\
   (0, m+1-k-l, l-m-1)  & {\text{if }} (k,l,m)\in {\rm (iii)},  \\
   (0, 0, -k)  & {\text{if }} (k,l,m)\in {\rm {(iv)}},
  \end{array} \right. 
\end{eqnarray*}
and 
\begin{align*}
{\rm (i)} &:\ \left\{(k,l,m); m \geq 1,\ m-k-l \leq -1 \right\},\ 
{\rm (ii)} :\ \left\{(k,l,m); m \geq 1,\ m-k-l \geq 0 \right\}, \\  
{\rm (iii)} &:\ \left\{(k,l,m); m \leq 0,\ m-k-l \leq -1 \right\},\ 
{\rm (iv)} :\ \left\{(k,l,m); m \leq 0,\ m-k-l \geq 0 \right\}  
\end{align*}
(cf. Propositions 3.4 and 3.9 in [Eb]).
Moreover, $q_0(x)$ and $r_0(x)$ can be expressed as sums of products of two hypergeometric series. 
For example, in the case of (i), $q_0(x)$ is given by
\begin{align}
q_0(x)&=q_1(a,b,c,k,l,m,x)\notag\\
&:=-\frac{(a,k)(b,l)}{(1-c)(c,m)}x^m{}_2F_1(c-a+m-k,c-b+m-l;c+m;x)\notag \\
&\times{}_2F_1(a+1-c, b+1-c;2-c;x)\notag \\
&+\frac{(a+1-c,k-m)(b+1-c,l-m)}{(1-c)(2-c,-m)}{}_2F_1(a, b; c;x)\notag\\
&\times{}_2F_1(1-a-k, 1-b-l; 2-c-m;x),
\end{align}
(cf. Theorem 3.7 in [Eb]). In addition, $q_0(x)$ and $r_0(x)$ possess certain symmetries. For example, in the case of (i),
$q_0(x)$ also has the expression
\begin{align}
q_0(x)&=x^{l-1}q_1(b+1-c,b,b+1-a,l-m,l,l-k,1/x)\notag\\
&=x^{l-1}\left[-\frac{(b+1-c,l-m)(b,l)}{(a-b)(b+1-a,l-k)}\left(\frac{1}{x}\right)^{l-k}\notag \right.\\
&\times{}_2F_1(c-a+m-k,1-a-k;b+1-a+l-k;1/x)\notag \\
&\times{}_2F_1(a+1-c,a;a+1-b;1/x)\notag\\
&+\frac{(a+1-c,k-m)(a,k)}{(a-b)(a+1-b,k-l)}{}_{2}F_{1}(b+1-c,b;b+1-a;1/x)\notag \\
&\left.\times _{2}F_1(c-b+m-l,1-b-l;a+1-b+k-l;1/x)\right]
\end{align}
(cf. Theorem 3.8 in [Eb]).
Now, we investigate $q_0(x)$ by comparing (4.6) with (4.7).
Focusing the denominators of (4.6) and (4.7), 
we find that the denominator of $q_0(x)$ has no term contained in
$$
\pm c+\Z ,\, \pm a\mp b+\Z .  
$$
Furthermore, noting that $q(x)$ can also be expressed as (4.5),
$2$, $3$, $4$, $\cdots$ are not contained as factors in the denominator of $q_0(x)$. 
Thus, for example, in the case of (i) with $k<0$,
we find that 
$$
q_0(x)=\frac{{\text{(a polynomial in $a,b,c,x$ over \Z)}}}{(1-a,-k)(c-a,m-k)}.
$$
We can obtain expressions for $q_0(x)$ in other cases similarly,
and we find the form of $r_0(x)$ is similar to that of $q_0(x)$:  
\begin{prp}
We assume $k\leq l$. 
Then, $q_0(x)$ and $r_0(x)$ can each be expressed as follows:
\begin{align*}
&\frac{{\text{$($a polynomial in $a,b,c,x$ over \Z$)$}}}{(1-a,-k)(c-a,m-k)}\quad
{\text for\ } {\rm (i)\ } {\text with\ } k<0, \\
&\frac{{\text{$($a polynomial in $a,b,c,x$ over \Z$)$}}}{(c-a,m-k)(c-b,m-l)}\quad
{\text for\ } {\rm (i)\ } {\text with\ } k\geq 0, k-m < 0, l-m<0,\\ 
&\frac{{\text{$($a polynomial in $a,b,c,x$ over \Z$)$}}}{(c-a,m-k)}\quad 
{\text for\ } {\rm (i)\ } {\text with\ } k\geq 0, k-m < 0, l-m\geq 0, \\ 
&{\text{$($a polynomial in $a,b,c,x$ over \Z$)$}}\quad
{\text for\ } {\rm (i)\ } {\text with\ } k\geq 0, k-m \geq 0, l-m\geq 0, \\
&\frac{{\text{$($a polynomial in $a,b,c,x$ over \Z$)$}}}{(1-a,-k)(1-b,-l)(c-a,m-k)(c-b,m-l)}\quad 
{\text for\ } {\rm (ii)\ } {\text with\ } l< 0, \\
&\frac{{\text{$($a polynomial in $a,b,c,x$ over \Z$)$}}}{(1-a,-k)(c-a,m-k)(c-b,m-l)}\quad 
{\text for\ } {\rm (ii)\ } {\text with\ } l\geq 0, k<0, l-m<0,\\
&\frac{{\text{$($a polynomial in $a,b,c,x$ over \Z$)$}}}{(1-a,-k)(c-a,m-k)}\quad 
{\text for\ } {\rm (ii)\ } {\text with\ } l\geq 0, k<0, l-m\geq 0,\\
&\frac{{\text{$($a polynomial in $a,b,c,x$ over \Z$)$}}}{(c-a,m-k)(c-b,m-l)}\quad 
{\text for\ } {\rm (ii)\ } {\text with\ } l\geq 0, k\geq 0, k-m<0, l-m< 0,\\
&\frac{{\text{$($a polynomial in $a,b,c,x$ over \Z$)$}}}{(c-a,m-k)}\quad 
{\text for\ } {\rm (ii)\ } {\text with\ } l\geq 0, k\geq 0, k-m<0, l-m\geq 0,\\
&\frac{{\text{$($a polynomial in $a,b,c,x$ over \Z$)$}}}{(1-a,-k)(1-b,-l)}\quad 
{\text for\ } {\rm (iii)\ } {\text with\ } k< 0, l< 0, k-m\geq 0,\\
&\frac{{\text{$($a polynomial in $a,b,c,x$ over \Z$)$}}}{(1-a,-k)(c-a,m-k)}\quad 
{\text for\ } {\rm (iii)\ } {\text with\ } k< 0, l\geq 0, k-m< 0,\\
&\frac{{\text{$($a polynomial in $a,b,c,x$ over \Z$)$}}}{(1-a,-k)}\quad 
{\text for\ } {\rm (iii)\ } {\text with\ } k< 0, l\geq 0, k-m\geq 0,\\
&{\text{$($a polynomial in $a,b,c,x$ over \Z$)$}}\quad 
{\text for\ } {\rm (iii)\ } {\text with\ } k\geq 0, l\geq 0, k-m\geq 0,\\
&\frac{{\text{$($a polynomial in $a,b,c,x$ over \Z$)$}}}{(1-a,-k)(1-b,-l)(c-a,m-k)(c-b,m-l)}\quad \\
&\hspace{7cm} {\text for\ } {\rm (ii)\ } {\text with\ } l<0, k-m<0, l-m<0, \\
&\frac{{\text{$($a polynomial in $a,b,c,x$ over \Z$)$}}}{(1-a,-k)(1-b,-l)(c-a,m-k)}\quad 
{\text for\ } {\rm (iv)\ } {\text with\ } l<0, k-m<0, l-m\geq 0, \\
&\frac{{\text{$($a polynomial in $a,b,c,x$ over \Z$)$}}}{(1-a,-k)(1-b,-l)}\quad 
{\text for\ } {\rm (iv)\ } {\text with\ } l<0, k-m\geq 0, l-m\geq 0, \\
&\frac{{\text{$($a polynomial in $a,b,c,x$ over \Z$)$}}}{(1-a,-k)(c-a,m-k)}\quad 
{\text for\ } {\rm (iv)\ } {\text with\ } l\geq 0, k-m<0, l-m\geq 0, \\
&\frac{{\text{$($a polynomial in $a,b,c,x$ over \Z$)$}}}{(1-a,-k)}\quad 
{\text for\ } {\rm (iv)\ } {\text with\ } l\geq 0, k-m\geq 0, l-m\geq 0.
\end{align*}
\end{prp}
Now, fix $(k,l,m)\in {\rm (iv)}$ with $l\geq 0, k-m<0, l-m\geq 0$.
Then, using (4.4), $Q(x)$ and $R(x)$ can be expressed as
\begin{gather*}
Q(x)=\frac{ab(c,m)x(1-x)Q_0'(x)}{c(b,l)(c-a,m-k)},\ 
R(x)=\frac{(c,m)R_0'(x)}{(b,l)(c-a,m-k)},
\end{gather*}
where $Q_0'(x), R_0'(x)\in \Z [a,b,c,x]$.
In particular, we consider the case that $(a,b,c)=(-2n_2,-2n_2,1), (k,l,m)=(-2n_1,2n_1,0)$.
Then, we obtain the following:
\begin{lmm}
Fix $n_1$. Then,
\begin{align}
\ell_1(n_1,n_2)&={}_2F_1(-2n_2-2n_1, -2n_2+2n_1; 1; -1)\notag \\
&=Q''{}_2F_1(-2n_2+1, -2n_2+1; 2; -1)+R''{}_2F_1(-2n_2, -2n_2; 1;-1),
\end{align}
where 
\begin{gather*}
Q''=\frac{-8n_2 ^2Q_0''}{(-2n_2,2n_1)(2n_2+1,2n_1)},\ 
R''=\frac{R_0''}{(-2n_2,2n_1)(2n_2+1,2n_1)},
\end{gather*}
and $Q_0'', R_0'' \in \Z[n_2]$. 
\end{lmm}

\subsection{Partial answer to Question 3}
In this subsection,
we investigate Question 3, using the result obtained in the previous subsection.

To begin, 
we obtain the closed form of $_2F_1(-2n_2+1, -2n_2+1; 2; -1)$ contained in the right hand side of (4.8).
Substituting $(a,b,c,x)=(-2n_2,-2n_2,1,-1)$ into the three term relation 
\begin{multline}
_2F_1 \left( a-1,b-1;c;x \right) =
{\frac {ab \left( c+1-a-b \right) 
 x\left( 1-x \right){}_2F_1 \left( a+1,b+1;c+1;x \right) }{ \left( c-a
 \right)  \left( c-b \right) c}}
 \\+{\frac { \left(  \left( {a}^{2}+{b}^{
2}- \left( a+b \right)  \left( c+1 \right) +ab+c \right) x+ \left( c-a
 \right)  \left( c-b \right)  \right) {}_2F_1 \left( a,b;c;x \right) }{
 \left( c-a \right)  \left( c-b \right) }},
\end{multline}
we obtain
\begin{multline}
_2F_1 \left( -2\,n_2-1,-2\,n_2-1;1;-1 \right) ={\frac {-16n_2^{2}{}_2F_1 \left( -2
\,n_2+1,-2\,n_2+1;2;-1 \right) }{2\,n_2+1}}\\
+{\frac {-4n_2\,{}_2F_1 \left( -2\,n_2,-2\,n_2
;1;-1 \right) }{2\,n_2+1}}.
\end{multline}
Note that from $(4.2)$, we know that the left hand side of (4.10) has a closed form.
Thus, $_2F_1(-2n_2+1, -2n_2+1; 2; -1)$ also has a closed form:
\begin{gather}
_2F_1(-2n_2+1, -2n_2+1; 2; -1)=\frac{(-1)^{n_2+1}\{1\cdot 3\cdot 5\cdots (2n_2-1)\}2^{n_2-2}}{(1,n_2)n_2}.
\end{gather}

Using Lemma 4.2 together with (4.3) and (4.11), we find
\begin{align}
&\ell_1(n_1,n_2)\notag \\
&=\frac{(-1)^{n_2}(2)^{n_2}
\{1\cdot 3\cdot 5\cdots (2n_2-1)\}(2n_2Q''_0+R''_0)}{(-2n_2,2n_1)(2n_2+1,2n_1)(1,n_2)}\notag\\
&=\frac{(-1)^{n_2}2^{n_2-2n_1}\{1\cdot 3\cdot 5\cdots (2n_2-2n_1-1)\}(2n_2Q''_0+R''_0)}{(1,n_2+n_1)(n_2+1-n_1,n_1)
\{(2n_2+1)(2n_2+3)(2n_2+5)\cdots (2n_2+2n_1-1)\}}.
\end{align}
Let us consider $\{1\cdot 3\cdot 5\cdots (2n_2-2n_1-1)\}$, appearing in the numerator of 
the right hand side of (4.12).
If there is a prime number $p$ satisfying $(n_2+n_1)<p< (2n_2-2n_1)$ in these factors,
then $\ell _1(n_1,n_2)$ is a multiple of $p$ (cf. Remark 2.1).
Thus, in this case, $\ell _1(n_1,n_2)$ is not a power of $2$,
and thus, $\ell _1(n_1,n_2)\neq r _1(n_1,n_2)$.
Therefore, we consider the case that $\pi(2n_2-2n_1)-\pi(n_2+n_1)>0$,
where the function $\pi (x)$ is defined as the number of primes up to $x$.
For example, we consider the case that $(n_1,n_2)=(n_1,4n_1+k)$,
where $k \in \Z _{\geq 0}$.
The following theorem provides an evaluation of $\pi(x)$:
\begin{thm}(cf. Section 4 in [Dus])
If $x\geq 599$, the following relation holds:
\begin{gather*}
\frac{x}{\log x}\left( 1+\frac{0.922}{\log x}\right)
\leq \pi (x)\leq
\frac{x}{\log x}\left( 1+\frac{1.2762}{\log x}\right).
\end{gather*}
\end{thm}
Let us define $x:=(5n_1+k)$,  
and suppose that $x\geq 599$. Then we find
\begin{align*}
&\pi(2n_2-2n_1)-\pi(n_2+n_1)=\pi (6n_1+2k)-\pi (5n_1+k)\\
\geq &\pi (1.2(5n_1+k))-\pi (5n_1+k)= \pi (1.2x)-\pi (x)\\
\geq &\frac{1.2x}{\log 1.2x}\left( 1+\frac{0.922}{\log 1.2x}\right)
-\frac{x}{\log x}\left( 1+\frac{1.2762}{\log x}\right)
>0
\end{align*} 
by applying Thorem 4.3.
Next, we assume $x< 599$. Then, we can determine whether $\ell _1 (n_1,n_2)$ is equal to $r_1 (n_1,n_2)$ 
using a computer-aided method, because there is only a finite number of pairs $(n_1,n_2)$ such that $x< 599$ is satisfied. 
Thus, we obtain a partial answer to Question 3:
There is no pair $(n_1,n_2)$ with $4n_1 \leq n_2$ satisfying (1.3) other than $(n_1,n_2)=(1,4)$.
Thus, we obtain the following:
\begin{prp}
For Question 1, 
there is no pair $(\kappa _1,\kappa _2)$ with $4\kappa _1 \leq \kappa _2$ satisfying (1.1)
other than $(\kappa _1,\kappa _2)=(4,16)$.
\end{prp}
\begin{rmk}
In this section, we considered the case that $\pi(2n_2-2n_1)-\pi(n_2+n_1)>0$.
We point out that $(n_1, n_2)$ satisfying this relation must satisfy 
$(2n_2 -2n_1)-(n_2+n_1)>0$, which implies $3n_1 < n_2$.
Thus, we cannot apply the method used in this section to the case $3n_1 \geq n_2$.
Therefore, we need other methods to answer Question 1 completely.
\end{rmk}
\textbf{Acknowledgement.}
The author would like to thank Professor Hiroaki Narita
for posing these interesting problems.
The author also thanks Dr. Genki Shibukawa
for many valuable comments.

\medskip
\begin{flushleft}
Akihito Ebisu\\
Department of Mathematics\\
Kyushu University\\
Nishi-ku, Fukuoka 819-0395\\
Japan\\
a-ebisu@math.kyushu-u.ac.jp
\end{flushleft}

\end{document}